

\baselineskip=14pt
\parskip=10pt
\def\halmos{\hbox{\vrule height0.15cm width0.01cm\vbox{\hrule height
  0.01cm width0.2cm \vskip0.15cm \hrule height 0.01cm width0.2cm}\vrule
  height0.15cm width 0.01cm}}

\magnification=\magstephalf

\def\1{{\overline{1}}}
\def\2{{\overline{2}}}
\parindent=0pt
\overfullrule=0in

\def\frac#1#2{{#1 \over #2}}
\centerline
{\bf An elegant Multi-Integral that implies an even more   }
\centerline
{\bf elegant determinant identity of Dougherty and McCammond  }
\bigskip
\centerline
{\it Tewodros AMDEBERHAN and Doron ZEILBERGER}

The first purpose of this note is to prove the following elegant identity.

{\bf Theorem A}: Let $z_1, \dots z_n$ be  commuting indeterminates, let $n$ be a positive integer, and let $a_1, \dots, a_n$, and $b$ be non-negative integers. Let
$\bar{a} := \sum_{i=1}^{\lceil n/2 \rceil} a_{2i-1}$. Then
$$
\int_0^{z_1} \dots  \int_0^{z_n} \, \prod_{i=1}^{n} x_i^b \prod_{1\leq j,k\leq n}(x_j -z_k)^{a_k} \prod_{1 \leq i<j \leq n} (x_j -x_i) \, dx_n \dots dx_1
$$
$$
\qquad =
\,(-1)^{\bar{a}} \prod_{1 \leq i<j \leq n} (z_j - z_i)^{a_i+a_j+1} \cdot \prod_{i=1}^{n} z_i^{a_i+b+1} \cdot
\frac{b! \, \prod_{i=1}^{n} a_i!}{(n+b+ \sum_{i=1}^{n} a_i)!}. \quad 
\eqno(1)
$$

The second purpose is to deduce from it (and thereby give a shorter proof) of the following even more elegant identity, discovered, and first proved in [1].

{\bf Theorem B (Dougherty and McCammond)}: Let
$$
p(Z) := \, \int_{0}^{Z} \, \prod_{i=1}^{n} (w-z_i)^{a_i} \, dw, \quad
$$
and let ${\bf J}(z_1, \dots, z_n)$ be the $n \times n$ matrix whose $(i,j)$ entry is
${\bf J}(z_1, \dots, z_n)_{i,j} \, := \, \frac{\partial}{\partial z_i} p(z_j),$
then
$$
\det \, {\bf J}(z_1, \dots, z_n) \, = \, 
\frac{\prod_{i=1}^{n} a_i!} { (\sum_{i=1}^{n} a_i)!}
\cdot \prod_{i=1}^{n} (-z_i)^{a_i}  \cdot \prod_{{{1 \leq i,j\leq n} \atop {i \neq j}}} (z_i - z_j)^{a_j}. \quad 
\eqno(2)
$$

{\bf Proof that ${\bf A} \Rightarrow {\bf B}$}: Let's rewrite the determinant and apply Cauchy's alternant formula so that
$$
\det ({\bf J}) 
= \det \left ( -a_i \int_0^{z_j} \prod_{k=1}^{n} (w-z_k)^{a_k} \frac{dw}{w-z_i} \right ) \, = \, \det \left ( -a_i \int_0^{z_j} \prod_{k=1}^{n} (x_j-z_k)^{a_k} \frac{d x_j}{x_j-z_i} \right ) $$
$$=\prod_{i=1}^{n} (-a_i) \int_0^{z_1} \dots \int_0^{z_n} \, \prod_{1 \leq j,k \leq n} (x_j -z_k)^{a_k} \cdot \det \left ( \frac{1}{x_j-z_i} \right)_{i,j}^{1,n} \, dx_n \cdots dx_1 $$
$$\qquad \qquad =\, \prod_{i=1}^{n} (-a_i) \int_0^{z_1} \dots \int_0^{z_n} \, \prod_{1 \leq j,k \leq n} (x_j -z_k)^{a_k} \cdot
\left [ \frac{ \prod_{1 \leq i < j \leq n} (z_i - z_j) (x_j -x_i)}{ \prod_{1 \leq i,j \leq n} (x_j -z_i)} \right ] \, dx_n \cdots dx_1$$
$$\qquad =\, \prod_{i=1}^{n} (-a_i)  \cdot \prod_{1 \leq i < j \leq n} (z_i - z_j)
\cdot \int_0^{z_1} \dots \int_0^{z_n}  \prod_{1 \leq j,k \leq n} (x_j -z_k)^{a_k-1}  \cdot \prod_{1 \leq i < j \leq n}  (x_j -x_i)  \, dx_n \cdots dx_1. \quad
$$
But by Theorem {\bf A}, with $b=0$ and $(a_1, \dots, a_n)$ replaced by $(a_1-1, \dots a_n-1)$, this equals
$$
=\, \,(-1)^{\bar{a}-\lceil n/2\rceil}
 \prod_{i=1}^{n} (-a_i)  \prod_{1 \leq i < j \leq n} (z_i - z_j)
\prod_{1 \leq i<j \leq n} (z_j - z_i)^{a_i+a_j-1}  \prod_{i=1}^{n} z_i^{a_i} 
\frac{ \prod_{i=1}^{n} (a_i-1)!}{(\sum_{i=1}^{n} a_i)!} \quad 
$$
$$
=\frac{\prod_{i=1}^{n} a_i!} { (\sum_{i=1}^{n} a_i)!}\,\,
\prod_{i=1}^{n} (-z_i)^{a_i}  \prod_{{{1 \leq i,j \leq n} \atop {i \neq j}}} (z_i - z_j)^{a_j}. \quad  \halmos 
$$

{\bf Proof of Theorem A}: The proof is by induction on $n$ and $b$. When $n=1$ and $b=0$ this is saying that $\int_0^{z_1} (x_1-z_1)^{a_1}\, dx_1 = (-1)^{a_1} \frac{(z_1)^{a_1+1}}{a_1+1}$, whose proof is left to
the reader's five-year-old.

Let's denote the statement of theorem {\bf A} by ${\bf A}(n,b)$.

{\bf Proof that ${\bf A}(n,b) \Rightarrow {\bf A}(n,b+1)$}

We claim that {\it both} sides of Eq. (1), let's call them $L(a_1, \dots, a_n;b)$ and  $R(a_1,\dots, a_n;b)$ respectively,
satisfy the recurrence
$$
X(a_1, \dots, a_n; b+1) \, = \,
\sum_{i=1}^{n} \left ( \prod_{{{j=1} \atop {j \neq i}}}^{n} \frac{z_j}{z_j-z_i} \right ) \cdot
X(a_1, \dots. a_i+1, \dots, a_n;b) 
+ \left ( \prod_{j=1}^{n} z_j \right) \cdot X(a_1, \dots, a_n; b). \quad
\eqno(3)
$$
In other words, if you replace $X$ by either $L$ or $R$ you get a true statement. Regarding the left side of $(1)$, in fact, this identity is already true if you replace $X$ by the {\bf integrand} of the left side of $(1)$,
since there are no $x_i$'s in sight, it is still true when you integrate with respect to $x_1, \dots, x_n$. We leave both checks as pleasant exercises for the reader. \halmos 

{\bf Proof that ${\bf A}(n-1,b)$ for all $b$ implies ${\bf A}(n,0)$ }

Fix $a_1, \dots, a_n$. Let $V(x_1, \dots, x_n) := \, \prod_{1 \leq i<j \leq n} (x_j -x_i)$ and
$$
F(x_1, \dots x_n; z_1, \dots , z_n):= \prod_{1\leq j,k\leq n} (x_j-z_k)^{a_k}. \quad
$$

We claim (check!) that
$$
\left(n+\sum_{i=1}^{n} a_i\right) \,F(x_1, \dots x_n; z_1, \dots, z_n) V(x_1, \dots,x_n)
$$
$$
=\sum_{i=1}^{n} (-1)^i V(x_1, \dots,x_{i-1}, x_{i+1}, \dots,  x_n) \cdot \frac{\partial}{\partial x_i} \left [ \prod_{j=1}^{n} (x_i-z_j) \cdot F(x_1, \dots , x_n; z_1, \dots, z_n) \right ]. \quad
\eqno(4)
$$

$$ $$

$$ $$

Applying  $\int_0^{z_1} \dots \int_0^{z_n} dx_n \cdots dx_1$, we get
$$
\left(n+\sum_{i=1}^{n} a_i\right)\,\int_0^{z_1} \dots \int_0^{z_n} F(x_1, \dots x_n; z_1, \dots, z_n) V(x_1, \dots,x_n) dx_n \cdots dx_1
$$
$$
= \sum_{i=1}^{n} (-1)^i \,\int_0^{z_1} \dots \int_0^{z_n}V(x_1, \dots,x_{i-1}, x_{i+1}, \dots,  x_n) \cdot \frac{\partial}{\partial x_i} \left [ \prod_{j=1}^{n} (x_i-z_j) \cdot F(x_1, \dots , x_n; z_1, \dots, z_n) \right ] dx_n \cdots dx_1 
$$
$$
= \, 
\sum_{i=1}^{n} (-1)^i \, 
\int_0^{z_1} \dots \int_0^{z_{i-1}} \int_0^{z_{i+1}} \dots \int_0^{z_n} V(x_1, \dots,x_{i-1}, x_{i+1}, \dots,  x_n) dx_n \cdots dx_{i+1} dx_{i-1} \cdots dx_1 
$$
$$ 
\times \int_{0}^{z_i} \frac{\partial}{\partial x_i} \left [ \prod_{j=1}^{n} (x_i-z_j) \cdot F(x_1, \dots , x_n; z_1, \dots, z_n) \right ] dx_i. \quad
$$
By the {\it Fundamental Theorem of Calculus}, we have
$$ 
\int_{0}^{z_i} \frac{\partial}{\partial x_i} \left [ \prod_{j=1}^{n} (x_i-z_j) \cdot F(x_1, \dots , x_n; z_1, \dots, z_n) \right ] dx_i
= \prod_{j=1}^{n} (x_i-z_j) \cdot F(x_1, \dots , x_n; z_1, \dots, z_n) \bigl \vert_{x_i=0}^{x_i=z_i}
$$
$$
= \,\prod_{j=1}^{n} (x_i-z_j)^{a_{j}+1}  \prod_{j=1}^{n}\prod_{{1 \leq i'\leq n} \atop {i' \neq i}} (x_i'-z_j)^{a_{j}} \bigl \vert_{x_i=0}^{x_i=z_i}
\,=\, - \prod_{j=1}^{n} (-z_j)^{a_{j}+1}  \prod_{j=1}^{n}\prod_{{1 \leq i'\leq n} \atop {i' \neq i}} (x_i'-z_j)^{a_{j}}.  \quad
$$
Going back we have that the left side of Eq. (1), when $b=0$,  is
$$
\frac{-1}{n+\sum_{i=1}^{n} a_i} \cdot
\prod_{j=1}^{n} (-z_j)^{a_{j}+1}  \cdot
\, \sum_{i=1}^{n} (-1)^i \,\int_0^{z_1} \dots  \int_{0}^{z_{i-1}}  \int_{0}^{z_{i+1}}  \dots \int_0^{z_n}
V(x_1, \dots,x_{i-1}, x_{i+1}, \dots,  x_n) 
$$
$$
\times \prod_{j=1}^{n}\prod_{{1 \leq i'\leq n} \atop {i' \neq i}} (x_i'-z_j)^{a_{j}} dx_n \cdots dx_{i+1} dx_{i-1} \cdots dx_1 
$$
$$
= \, \frac{-1}{n+\sum_{i=1}^{n} a_i} \cdot \prod_{j=1}^{n} (-z_j)^{a_{j}+1}  \cdot
\, \sum_{i=1}^{n} (-1)^i \,\int_0^{z_1} \dots  \int_{0}^{z_{i-1}}  \int_{0}^{z_{i+1}} \dots \int_0^{z_n}V(y_1, \dots, y_{n-1}) 
$$
$$
\times \prod_{j=1}^{n}\prod_{i=1}^{n-1} (y_i-z_j)^{a_{j}} dy_{n-1} \cdots  dy_1. \quad
$$

$$ $$

We now claim that
$$
\, \sum_{i=1}^{n} (-1)^i \,\int_0^{z_1} \dots  \int_{0}^{z_{i-1}}  \int_{0}^{z_{i+1}} \dots \int_0^{z_n}
V(y_1, \dots, y_{n-1}) \cdot
\prod_{j=1}^{n}\prod_{i=1}^{n-1} (y_i-z_j)^{a_{j}} dy_{n-1} \cdots  dy_1
$$
$$
=\,\int_{z_1}^{z_2} \dots \int_{z_1}^{z_n} V(y_1, \dots , y_{n-1}) \cdot \prod_{j=1}^{n}\prod_{i=1}^{n-1} (y_i-z_j)^{a_{j}} dy_{n-1} \cdots  dy_1. \quad 
\eqno(5)
$$

In order to prove this, notice that each of the integrands on the left, and the integrand on the right, are {\it anti-symmetric} in their arguments. Hence, for any given permutation
of the integration variables, the effect is to multiple it by the sign of that permutation. Calling the common integrand $f(y_1,\dots, y_{n-1})$ and denoting $A_n(i)=Per(1,\dots, i-1, i+1, \dots n)$, we claim that
$$
\, \sum_{i=1}^{n} (-1)^i \,
\sum_{\pi \in A_n(i)} \, sgn(\pi) \, \int_0^{z_{\pi(1)}} \dots  \int_{0}^{z_{\pi(i-1)}} \int_{0}^{z_{\pi(i+1)}} \cdots \int_0^{z_{\pi(n)}}
f(y_1, \dots, y_{n-1}) dy_{n-1}\cdots dy_1
$$
$$
=\,\sum_{\pi \in A_n(1)}  sgn(\pi) \, \int_{z_1}^{z_{\pi(2)}} \dots \int_{z_1}^{z_{\pi(n)}} f(y_1, \dots, y_{n-1}) dy_{n-1} \cdots  dy_1. \quad 
\eqno(6)
$$
Since both sides of Eq. $(6)$ are $(n-1)!$ times the respective sides of Eq. $(5)$, if we can prove $(6)$, then $(5)$ would follow.

But {\bf surprise!}, Eq. $(6)$ is valid for {\it any} integrand! It is just a relation between {\it regions} in $R^{n-1}$ that is equivalent to an easy symmetric function identity,
that we also leave as a pleasant exercise to the reader. Now make the change of variables $(y_1, \dots, y_{n-1}) \rightarrow (y_1-z_1, \dots, y_{n-1}-z_1) $, thereby making
it a case of ${\bf A}(n-1,b)$ with $b=a_1$;  and $a_1, \dots, a_{n-1}$ replaced by $a_2, \dots, a_n$, respectively; and
 $z_1, \dots, z_{n-1}$ replaced by $z_2-z_1, \dots, z_n-z_1$, respectively. Plugging it in and simplifying, completes
the induction. $\halmos $ $\halmos $

{\bf Comment}: Readers that prefer not do the `exercises' can convince themselves of all the claims, empirically, by playing with the Maple package {\tt CritVal.txt} available
from the front of this paper

{\tt https://sites.math.rutgers.edu/\~{}zeilberg/mamarim/mamarimhtml/crit.html} \quad.

\bigskip
{\bf Reference}

[1] M. Dougherty, J. M. McCammond, {\it Critical Points, critical values, and a determinant identity for complex polynomials}, Proc. Amer. Math. Soc. {\bf 148} (2020), 5722-5289. \hfill\break
{\tt https://arxiv.org/abs/1908.10477} \quad .

\bigskip

\hrule

$$ $$

\bigskip
Tewodros Amdeberhan, Department of Mathematics, Tulane University, New Orleans, LA 70118
Email: {\tt tamdeber@tulane.edu}   \quad

Doron Zeilberger, Department of Mathematics, Rutgers University (New Brunswick), Hill Center-Busch Campus, 110 Frelinghuysen
Rd., Piscataway, NJ 08854-8019, USA. \hfill\break
Email: {\tt DoronZeil@gmail.com}   \quad 

\bigskip
\hrule
\bigskip
Jan. 7, 2022

\end